\documentclass[12pt]{article}
\usepackage{amssymb}
\usepackage{a4wide}
\usepackage{epsfig}
\usepackage{oldgerm}
\usepackage{amsmath}
\usepackage{cite}
\usepackage{amsthm}
\usepackage{fixmath}
\usepackage{multicol}
\usepackage{amstext,color}
\usepackage{array,multirow}
\usepackage{fullpage}
\usepackage{graphicx,graphics}
\usepackage{caption}
\usepackage{subcaption}

\usepackage[vlined,ruled,linesnumbered]{algorithm2e}

\definecolor{teal}{rgb}{0.0,0.5,0.5}
\definecolor{frenchrose}{rgb}{0.96,0.29,0.54}
\definecolor{lg}{rgb}{0.36,0.99,0.82}
\definecolor{dblue}{rgb}{0,0,.5}
\definecolor{dpink}{cmyk}{.2,1,.1,.04}
\definecolor{purple}{rgb}{0.35,0.04,0.64}
\definecolor{borange}{rgb}{1, .388, 0}
\definecolor{dpurple}{rgb}{0.61,0.22,1.00}
\definecolor{purp}{rgb}{0.44,0.00,0.87}
\definecolor{green}{rgb}{0.00,0.44,0.00}
\definecolor{junebud}{rgb}{0.74,0.85,0.34}
\definecolor{plum}{rgb}{0.56,0.27,0.52}
\theoremstyle{plain}
\newtheorem{theorem}{Theorem}[section]
\newtheorem{corollary}[theorem]{Corollary}

\newtheorem{definition}[theorem]{Definition}
\newtheorem{lemma}[theorem]{Lemma}
\newtheorem{remark}[theorem]{Remark}

\newtheorem{prop}[theorem]{Proposition}
\newtheorem{example}[theorem]{Example}
\def\CC{{\textmd \kern.24em \vrule width.02em height1.4ex depth-.05ex \kern-.26emC}}
\def\TagOnRight
\def\QQ{\rlap {\raise 0.4ex \hbox{$\scriptscriptstyle |$}} {\hskip -0.1em Q}}
\catcode`\@=11
\begin{document}

\begin{center}
  {{\bf \large {\rm {\bf $\alpha$-Robust Error Analysis of $L2$-$1_{\sigma}$ Scheme on Graded Mesh for Time-fractional Nonlocal Diffusion Equation }}}}
\end{center}
\begin{center}
	{\textmd {\bf Pari J. Kundaliya}}\footnote{\it Department of Mathematics,  Institute of Infrastructure, Technology, Research And Management, Ahmedabad, Gujarat, India, (pariben.kundaliya.pm@iitram.ac.in)}
\end{center}

\begin{abstract}
In this work, a time-fractional nonlocal diffusion equation is considered. Based on the $L2$-$1_{\sigma}$ scheme on a graded mesh in time  and the standard finite element method (FEM) in space, the fully-discrete $L2$-$1_{\sigma}$ finite element method is investigated for a time-fractional nonlocal diffusion problem. We prove the existence and uniqueness of fully-discrete solution. The $\alpha$-robust error bounds are derived, i.e. bounds remain valid as $\alpha$ $\rightarrow {1}^{-},$ where $\alpha \ \in (0,1)$ is the order of a temporal fractional derivative. The numerical experiments are presented to justify the theoretical findings.
\end{abstract}


\noindent {\bf AMS(MOS):} 65M60, 65M15, 35R11.

\section{Introduction}

  \indent
  In recent years, many researchers have focused on standard nonlocal parabolic partial differential equations (PDEs) (see \cite{[me],[sk2017],[rbl1],[r6new],[zc]}). In this work, we consider the following time-fractional nonlocal diffusion equation:
   \begin{subequations}\label{2e1}
  	\begin{align}
  	\label{cuc3:1.1}
    \qquad  ^{c}_{0}{D}^{\alpha}_{t}u(x,t)-a(l(u)) \: \Delta u(x,t) & = f(x,t)  & \mbox{in }  & \quad \Omega\times(0,T], \\
  	\label{cuc3:1.2}
   	u(x,t) & = 0  & \mbox{on }  & \quad \partial\Omega\times(0,T],\\
   	\label{cuc3:1.3}
  	u(x,0)&=u_0(x) & \mbox{in } & \quad \Omega,
  	\end{align}
  \end{subequations}
  where $\Omega \subseteq \mathbb{R}^d$ $(d=1$ or $2)$ is a bounded domain with smooth boundary $\partial \Omega$, the function $a : \mathbb{R} \longrightarrow \mathbb{R}$ is a diffusion coefficient, and $l(u)= \int_{\Omega} u(x,t) \, dx.$ The term $^{c}_{0}{D}^{\alpha}_{t}u(x,t)$ is the $\alpha ^{th}$ order Caputo fractional derivative of $u$ with $\alpha \in (0,1)$ and it is defined in \cite[Definition 3.1]{[r1]} as \\
  \begin{equation*}\label{2e2}
  \begin{split}
  {^{c}_{0}D^{\alpha}_{t}}u(x,t) =& \frac{1}{\Gamma (1-\alpha)}\int_{0}^{t}(t-s)^{-\alpha} \; \frac{\partial u(x,s)}{\partial s} \, ds , \;\; t>0.
  \end{split}
  \end{equation*}
  The problem considered here is nonlocal due to the presence of the term $a(l(u)),$ which is determined by a global quantity. Nonlocal problems often lead to a number of mathematical and computational difficulties. This type of nonlocal problem appears in population dynamics, where the diffusion coefficient $a$ depends on the entire population in the domain rather than on a local density \cite{[zc],[sk2017],[sp1]}. \\
  \indent
  The solution of time-fractional PDEs in general exhibits weak singularity near time $t=0$ \cite{[kwn],[r02],[r08],[r8],[r13],[cmb]}. This makes the existing numerical methods with uniform time mesh often less accurate. In \cite{[r8]}, the authors consider a semilinear time-fractional diffusion equation and show that if the solution has singularity near time $t=0$, then the $L1$ scheme on uniform mesh (with step size $\tau$) attains $\mathcal{O}(\tau^\alpha)$ order of convergence in $L^{\infty}$-norm in time. The graded temporal meshes offer an efficient way of computing a reliable numerical solution near time $t=0$ \cite{[kwn],[r02],[r11]}. The $L1$ scheme is considered in \cite{[sp1]} on graded mesh for solving the time-fractional nonlocal diffusion problem. It is shown that, in the case of weak singularity of solution $u$ near the initial time $t = 0$, the optimal convergence result with ${\mathcal{O}}{(\tau^{2-\alpha})}$ in temporal direction is achieved. It is also known that $L1$ scheme can give at most $2-\alpha$  order of convergence in temporal direction. Therefore, in the present work, we investigate a higher order scheme; namely, the $L2$-$1_{\sigma}$ scheme. This scheme was first proposed by Alikhanov \cite{[AAl2]} on uniform mesh and latter this scheme was applied in~\cite{[hr12],[r16],[fmr13],[hwj1]} on general meshes. In \cite{[sp1]}, the error estimates are established, which are $\alpha$-nonrobust although the numerical scheme performs well as $\alpha$ $\rightarrow {1}^{-}.$  In present work, the $\alpha$-robust error bounds are derived in $L^2$-norm and $H^1$-norm. The concept of $\alpha$-robust bound was first introduced in \cite{[r15]} and then it was used in \cite{[cmb],[r16],[hwh]} for different problems.\\
  \indent
  Notation: Henceforth, $C > 0$ denotes a generic constant, which is independent of mesh parameters $h$ and $N.$  We write $(\cdot,\cdot)$ for the inner product and $ \|\cdot\|$ for the norm on $L^2(\Omega).$  For $m \in \mathbb{N}$, the notation $H^m(\Omega)$ denotes the standard Sobolev space with the norm $\|\cdot\|_m$ and the space $H_{0}^1(\Omega)$ is consist of functions from $H^1(\Omega)$ whose trace vanish on the boundary $\partial\Omega.$ \\
  \indent
  Throughout, the following hypotheses are made on the given data $a$, $f$ and $u_0$:
\begin{enumerate}
	\item[H1:] There exist $m_1, m_2 >0$ such that  $0<m_1 \le a(x) \le m_2 < \infty, \; \forall \, x \in \mathbb{R}.$
	\item[H2:] $a$ is Lipschitz continuous, i.e., $|a(x_1) - a(x_2)| \le L  |x_1-x_2|, \; (x_1, x_2 \in \mathbb{R})$.
	\item[H3:] $f \in L^{2}(0, T; L^2(\Omega))$ and $u_0 \in H_0^1(\Omega)\cap H^2(\Omega)$.
\end{enumerate}

  \indent
In Section 2, we give a fully-discrete scheme for problem \eqref{2e1}. We derive \emph{a priori} bound for a fully-discrete solution in Section 3. The existence and uniqueness of fully-discrete solution are proved in Section 4. Error analysis of our proposed scheme is presented in Section 5. Finally, our theoretical results are confirmed in Section 6.
  

\section{Fully-discrete Scheme}
\indent
In this section, we discretize the differential equation \eqref{2e1} in both space and time. Firstly, we write the weak formulation of problem \eqref{2e1} and it is given below.\\
For each $t \in(0,T],$ find $u(\cdot,t) \in  H^1_0(\Omega)$ such that
\begin{subequations}\label{2e3}
\begin{align}
 \label{2e3.a}
   (^{c}_{0}{D}^{\alpha}_{t}u, v ) \, + \, a\big(l(u)\big) \, (\nabla u, \nabla v) \, &=& \! \! \! \! \! \! \! \! \! \! \! \! \! \! \! \! \! \! \! \! \! \big(f, v\big), &\quad  \forall v \in  H^1_0(\Omega).& \\
 \label{2e3.b}
   u(\cdot, 0) &=& \! \! \! \! \! \! \! \! \! \! \! \! \! \! \! \! \! \! \! \! u_0 (\cdot), \: &\quad \mbox{in} \; \, \Omega. \quad&
\end{align}
\end{subequations}
Note that the above problem \eqref{2e3} has a unique weak solution under the Hypotheses H1-H3 stated in the introduction~\cite[Theorems 2.4 and 2.5]{[mn3]}.\\
\indent
In order to establish a fully-discrete scheme, we assume that $\tau_h$ is a quasi uniform partition of $\Omega$ into disjoint intervals in $\mathbb{R}^1$ or triangles in $\mathbb{R}^2$ and $h$ is the mesh size. Let $V_h$ be the $M$-dimensional subspace of $H^{1}_{0}(\Omega)$, which consists of continuous piecewise linear polynomials.\\
\indent
Let $N \in \mathbb{N}$ and $\tau_N := \Big\{ t_n : t_n = T\Big(\frac{n}{N}\Big)^{r}, \, n=0, \ldots, N\Big\}$
be a non-uniform partition of $[0,T]$ with time step $\tau_n = t_n - t_{n-1}$ and mesh grading parameter $r \ge 1$. Let $u^n := u(t_n)$ and $U^n$ be the approximate value of $u$ at $t_n$. For $\sigma = \frac{\alpha}{2}$, we set $t_{n-\sigma} := (1- \sigma) \, t_n + \sigma \, t_{n-1}$ and  $U^{n, \sigma} :=  (1- \sigma) \, U^{n} + \sigma \, U^{n-1} .$\\
\indent The $L2$-$1_{\sigma}$ approximation to $^c_0D^{\alpha}_{ t_{n-\sigma}}u$ on a graded mesh is given below~\cite[P.5]{[r16]}. For $n = 1, \ldots, N,$
\begin{eqnarray}\label{2e39}
^c_0D^{\alpha}_{ t_{n-\sigma}}u & \approx &  D^{\alpha}_{N} u^{n-\sigma}\nonumber\\
& := & \sum_{j=1}^{n} g_{n,j} (u^{j} - u^{j-1}) \nonumber\\
& = &  g_{n,n} \, u^{n} - g_{n,1} \, u^{0} - \sum_{j=2}^{n} \, \big(g_{n,j} - g_{n,j-1}\big) \, u^{j-1}.
\end{eqnarray}
Here $g_{1,1} = \tau_1^{-1} \, a_{1,1}$ and for $n \ge 2,$
\begin{eqnarray}\label{2e40}
g_{n,j} = \left\{
\begin{array}{lll}
\tau_{j}^{-1} \, (a_{n,1}-b_{n,1}) & (\text{if } j=1) \\
\tau_{j}^{-1} \, (a_{n,j}+b_{n,j-1}-b_{n,j}) & (\text{if } \; 2 \le j \le n-1) \\
\tau_{j}^{-1} \, (a_{n,n}+b_{n,n-1}) & (\text{if } j=n),
\end{array}
\right.
\end{eqnarray}
where
\begin{eqnarray}
	a_{n,n} & = & \frac{1}{\Gamma{(1-\alpha)}} \, \int_{t_{n-1}}^{t_{n-\sigma}} \, (t_{n-\sigma} - \eta)^{- \alpha} \, d\eta,  \; \; \mbox{for} \; n \ge 1. \nonumber \\
	a_{n,j} & = & \frac{1}{\Gamma{(1-\alpha)}} \, \int_{t_{j-1}}^{t_j} \, (t_{n-\sigma} - \eta)^{- \alpha} \, d\eta, \; \; \mbox{for}  \; 1 \le j \le n-1. \nonumber \\
	b_{n,j} & = & \frac{2}{\Gamma{(1-\alpha)} \,(t_{j+1} - t_{j-1})} \, \int_{t_{j-1}}^{t_j} \, (t_{n-\sigma} - \eta)^{- \alpha} \, (\eta - t_{j-\frac{1}{2}}) \, d\eta, \, \mbox{ for }  \; 1 \le j \le n-1. \nonumber
\end{eqnarray}

\begin{lemma}\label{2l6}
   Let $v(t) \in C[0,T] \cap C^3(0,T]$ be a function with $|v^{(q)}(t)| \le C (1 + t^{\alpha - q}),$ for $q=0,1,2,3$ and for $t \in (0,T]$. Then, for $n=1, \ldots, N,$
	\begin{equation*}\label{2e110}
	  t_{n-\sigma}^{\alpha} \, \| ^c_0D^{\alpha}_{ t_{n-\sigma}}v - D^{\alpha}_{N} \, v^{n-\sigma} \| \le C \, N^{-\min \{3- \alpha, \, r \alpha\}}.
	\end{equation*}		
\end{lemma}
\textbf{Proof.} The proof follows from \cite[Lemmas 1 and 7]{[fmr13]}.\\

\indent
Using the above notations, the \textit{fully-discrete scheme} for \eqref{2e1} is following:\\ 
For each $n=1, \ldots, N$, find $U^{n} \in V_h$ such that $\forall v_h \in V_h,$
\begin{subequations}\label{2e6}
	\begin{align}
  \label{2e6.a}
	\left( {D}^{\alpha}_{N}U^{n-\sigma}, v_h \right) \, + \, a\big(l(U^{n, \sigma})\big) \, (\nabla U^{n, \sigma}, \nabla v_h) =& \, (f^{n-\sigma}, v_h),\\
  \label{2e6.b}
	U^0 =& \, u^0_h,
	\end{align}
\end{subequations}
where $u^0_h$ is some approximation of $u_0(x)$.\\

Now, from the definition of ${D}^{\alpha}_{N}$ and as $\phi_i \,^,$s are the basis functions for $V_h$, from \eqref{2e6.a}, for each $i= 1,\ldots, M,$ we have
\begin{eqnarray}\label{2e47}
g_{n,n}( U^{n}, \, \phi_i) - \Big(g_{n,1} U^0 + \sum_{k=2}^{n} (g_{n,k} - g_{n,k-1}) U^{k-1}, \, \phi_i \Big) \!\!\! & + & \!\!\! a\big(l(U^{n, \sigma})\big)
(\nabla U^{n, \sigma}, \, \nabla \phi_i)  \nonumber \\
& = &\big(f^{n-\sigma}, \, \phi_i \big).
\end{eqnarray}

Since $U^{n} \in V_h,$ there exists $\alpha^{n}_j \in \mathbb{R}$ (for each $j = 1, \ldots, M$) such that
\begin{equation*}\label{2e48}
U^{n} \,= \, \sum_{j=1}^{M}  \alpha^{n}_j \phi_j.
\end{equation*}

Set $\alpha^{n} \, = [ \alpha^{n}_1, \alpha^{n}_2,\ldots,\, \alpha^{n}_M].$ It can be noticed that \eqref{2e47} is a nonlinear equation in $U^{n}.$ If we directly use the Newton's method for solving the equation \eqref{2e47}, then the Jacobian matrix becomes dense due to the presence of nonlocal term $a(l(u))$ \cite{[sp1]}. Due to this term, the scheme becomes computationally expensive. To overcome this difficulty, we follow the idea of \cite{[r3]}. This technique was first proposed for elliptic nonlocal problems and latter it was analyzed in \cite{[sk],[me],[sp1]} for solving nonlocal parabolic problems. For this, we reformulate \eqref{2e47} as follows: Find $U^{n}  \in V_h$ and $d \in \mathbb{R}$ such that
\begin{equation}\label{2e52}
  l(U^{n, \sigma}) - d \, = \, 0,
\end{equation}
and for $1 \le i \le M,$
\begin{eqnarray}\label{2e53}
(U^{n}, \, \phi_i) & - & \frac{1}{g_{n,n}} \, \Big(g_{n,1} U^0 + \sum_{k=2}^{n} (g_{n,k} - g_{n,k-1}) U^{k-1}, \, \phi_i \Big) \nonumber\\
& + & \frac{1}{g_{n,n}} \, a(d) (\nabla U^{n, \sigma}, \, \nabla \phi_i)  - \frac{1}{g_{n,n}} \, \big(f^{n-\sigma}, \, \phi_i\big) = 0.
\end{eqnarray}

Now, in order to use the Newton's method, we rewrite \eqref{2e52} and \eqref{2e53} as follows:
\begin{equation}\label{2e54}
  F_i(U^n, \, d) \, = \, 0, \quad \mbox{for} \; 1 \le i \le M+1,
\end{equation}	
where, for each $1 \le i \le M$,
\begin{eqnarray}\label{2e56}
  F_i(U^{n}, d) & = & \! ( U^{n}, \phi_i) - \frac{1}{g_{n,n}} \Big(g_{n,1} U^0 + \sum_{k=2}^{n} (g_{n,k} - g_{n,k-1}) U^{k-1},  \phi_i \Big)  \nonumber\\
  &  & + \frac{1}{g_{n,n}} \, a(d) (\nabla U^{n, \sigma}, \, \nabla \phi_i) - \frac{1}{g_{n,n}} \, \big(f^{n-\sigma}, \, \phi_i\big), \nonumber
\end{eqnarray}
and
\begin{equation*}\label{2e57}
 F_{M+1}(U^{n}, \, d) \, = \, l(U^{n, \sigma}) - d. 
\end{equation*}

An application of Newton's method in \eqref{2e54} gives the following system of equation:
\begin{gather*}\label{2e58}
 J
 \begin{bmatrix}
 {{\alpha}^{n}} \\ d
 \end{bmatrix}
 =
 \begin{bmatrix}
 A & b\\ c & -1
 \end{bmatrix}
 \begin{bmatrix}
 {{\alpha}^{n}} \\ d
 \end{bmatrix}
 =
 \begin{bmatrix}
 \bar{F} \\ F_{M+1}
 \end{bmatrix},
\end{gather*}
where $J$ is the Jacobian matrix, ${{\alpha}^{n}} \, = \, [ \alpha^{n}_1, \ldots,\, \alpha^{n}_M]^T,\,$  $\bar{F} \, = \, [F_1, \ldots, F_M ]^T$, and the entries of the matrices $ A_{M \times M},$ $b_{M \times 1}$ and $c_{1 \times M}$ are given below. 
\begin{eqnarray}\label{2e59}
  A_{ij}  & = &  ( \phi_j, \, \phi_i) + \frac{\sigma}{g_{n,n}} \, a(d) (\nabla \phi_j, \, \nabla \phi_i),  \nonumber \\
  b_{i1}  & = &  \frac{1}{g_{n,n}} \, a'(d) \, (\nabla U^{n, \sigma}, \, \nabla \phi_i), \nonumber \\
  c_{1j}  & = &  \sigma \, \int_{\Omega} \phi_j \, dx. \nonumber
\end{eqnarray}

Now, the sparsity of $J$ follows from the sparsity of matrix $A$. Note that the solutions of \eqref{2e47} and \eqref{2e52}-\eqref{2e53} are equivalent \cite[Theorem 3.1]{[sp1]}.

\section{\emph{A priori} Bound for Fully-discrete Solution}

\indent
In this section, we provide \emph{a priori} bound for fully-discrete solution of problem \eqref{2e6}. In order to this, we shall need the definition of the Discrete Laplacian and the coercivity property of $L2$-$1_{\sigma}$ scheme.

\begin{definition} \cite[P.10]{[vth]}
	The discrete \textit{Laplacian} is a map $\Delta _h : \, V_h \, \rightarrow \, V_h $ defined as
	\begin{equation*}\label{2e106}
	(\Delta _hu, \, v) \, = \, -(\nabla u, \, \nabla v), \quad \forall \: u, \, v \, \in \, V_h.
	\end{equation*}
\end{definition}

\begin{lemma}\label{2l1}
	\cite[Corollary 1]{[fmr13]} For any function $v = v( \cdot , \, t)$ defined on a graded mesh $\left\lbrace  t_j \right\rbrace ^{N}_{j=0}$, one has
	\begin{equation*}\label{2ee12}
	   \frac{1}{2} \, {D}^{\alpha}_{N} \|v^{n - \sigma}\|^2 \, \le \, \left( {D}^{\alpha}_{N} v^{n-\sigma}, \, v^{n, \sigma} \right), \quad \mbox{for} \; \; n=1, \ldots, N.
	\end{equation*}
\end{lemma}

For further analysis, for $n=1,\ldots,N$, we define the coefficients $p^{(n)}_{n-i}$ as follows:
\begin{eqnarray}\label{2e121}
p^{(n)}_{n-i} := \left\{
\begin{array}{ll}
\frac{1}{g_{i,i}} \sum_{k=i+1}^{n} \left( g_{k, i+1} - g_{k, i} \right) p^{(n)}_{n-k} & (\text{if } i = 1, \ldots , n-1)\\
\frac{1}{g_{n,n}} & (\text{if } \; i = n).
\end{array}
\right.
\end{eqnarray}

\begin{lemma}\label{2l7}
	\cite[P.46 (15)]{[r16]} Let $\gamma \in (0,1)$ be any constant. For $n=1,\ldots,N$, the coefficients defined in \eqref{2e121} satisfy the following inequality.
	\begin{equation*}\label{2e124}
	\begin{split}
	\sum_{i=1}^{n} \, p^{(n)}_{n-i} \, i^{r (\gamma - \alpha)} \, \le \, \frac{11 \, \Gamma(1 + \gamma - \alpha)}{4 \, \Gamma(1+ \gamma)} \, T^{\alpha} \, \bigg(\frac{t_n}{T}\bigg)^{\gamma} \, N^{r (\gamma - \alpha)}.\\
	\end{split}
	\end{equation*}	
\end{lemma} 	

Now, we state the discrete fractional Gr$\ddot{{o}}$nwall inequality, which will be used in deriving \emph{a priori} bounds and error estimates.
\begin{lemma}\label{2l5}
\cite[ Lemma 4.1]{[r16]} Let $n \ge 1$ and let $\lambda_i$ be non-negative constants with $0 < \sum_{i=0}^{n} \lambda_i \le \Lambda$, where $\Lambda$ is some positive constant independent of $n$. Suppose that the non-negative sequences $\{\xi ^n, \: \eta ^n : n \ge 1\}$ are bounded and the non-negative grid function $\{v^n : n \ge 0\}$ satisfies
 \begin{equation*}\label{2e105}
\begin{split}
D_N^{\alpha} (v^{n-\sigma})^2 \, \le \, \sum_{i=1}^{n} \, \lambda_{n-i} \, (v^{i, \sigma})^2 + \, \xi^{n} \, v^{n, \sigma} \, + \,  (\eta ^{n})^2, \quad \mbox{for} \; \, n \ge 1.
\end{split}
\end{equation*}
If the non-uniform grid satisfies the criterion $\max\limits_{1 \le n \le N}\tau_n \le \frac{1}{\sqrt[\alpha]{\frac{11}{2} \, \Gamma (2- \alpha) \, \Lambda}},$ then
\begin{equation*}\label{2e109}
 \begin{split}
    v^{n} \le 2 E_{\alpha} \bigg(\frac{11}{2} \, \Lambda \, t_{n}^{\alpha}\bigg) \, &\Big[ v^0 + \max_{1 \le k \le n} \sum_{j=1}^{k} \, p^{(k)}_{k-j} (\xi^j + \eta^j) + \max_{1 \le j \le n} \, \left\lbrace \eta^j \right\rbrace \Big].\\
 \end{split}
\end{equation*}	
\end{lemma}

\begin{remark}
   The $L2$-$1_{\sigma}$ scheme on graded mesh satisfies all the hypotheses, which are required for the proof of Lemmas~ \ref{2l7} and \ref{2l5} \cite[Example 3, P.222]{[r5]}.
\end{remark} 

The next theorem is one of the main results of this article.
\begin{theorem}\label{2th4}
	The fully-discrete solution $U^{n}$ of \eqref{2e6} satisfies
	\begin{eqnarray}
	\label{2e101}
   	  \|U^{n}\| & \le & C \, \big(1+\|U^0\| \big), \\
    \label{2e102}
	  \|\nabla U^{n}\| & \le & C  \, \big(1 + \|\nabla U^0\| \big),
	\end{eqnarray}
for each  $n=1, \ldots, N$.	
\end{theorem}

\noindent	
\textbf{Proof.} First we note that, by substituting $\gamma = \alpha$ in Lemma~\ref{2l7}, we can get
\begin{equation}\label{2e14a}
\begin{split}
  \sum_{j=1}^{k} p^{(k)}_{k-j} \, \le \,  \frac{11 \, t_k^{\alpha}}{4 \, \Gamma{(1+\alpha)}} \, \le \, \frac{11 \, T^{\alpha}}{4 \, \Gamma{(1+\alpha)}}.
\end{split}
\end{equation}
Next, by taking $v_h= U^{n, \sigma}$ in \eqref{2e6.a}, we get
\begin{equation}\label{2e8}
\begin{split}
  \left( {D}^{\alpha}_{N}U^{n-\sigma}, U^{n, \sigma} \right) \, + \, a\big(l(U^{n, \sigma})\big) \, (\nabla U^{n, \sigma}, \nabla U^{n, \sigma}) \, =& \, (f^{n-\sigma}, U^{n, \sigma}), \; \,  \forall v_h \in V_h. \\
\end{split}
\end{equation}

\noindent
Applying Hypothesis H1 and the Cauchy-Schwartz inequality in \eqref{2e8}, we have
\begin{equation}\label{2e9}
\begin{split}
  \left( {D}^{\alpha}_{N}U^{n-\sigma}, U^{n, \sigma} \right) \, \le& \, \|f^{n-\sigma}\| \, \|U^{n, \sigma}\|.
\end{split}
\end{equation}

\noindent
Using Lemma \ref{2l1} and the inequality $b_1b_2 \le \frac{b_1^2}{2} + \frac{b_2^2}{2}$ in \eqref{2e9}, we get
\begin{equation}\label{2e13}
\frac{1}{2} \, {D}^{\alpha}_{N} \|U^{n - \sigma}\|^2 \, \le \, \frac{1}{2} \, \|f^{n-\sigma}\|^2 + \frac{1}{2} \, \|U^{n, \sigma}\|^2.
\end{equation}

\noindent
An application of Lemma \ref{2l5} into \eqref{2e13} gives	
\begin{equation*}\label{2e14}
\begin{split}
\|U^{n}\|  \le& 2  E_{\alpha}\Big(\frac{11}{4} t_{n}^{\alpha}\Big) \Big( \|U^0\|  +  \! \max_{1 \le k \le n} \! \Big( \max_{1 \le j \le k} \|f^{j-\sigma}\|\Big) \sum_{j=1}^{k} p^{(k)}_{k-j} + \! \max_{1 \le j \le n}  \|f^{j-\sigma}\| \Big).\\
\end{split}
\end{equation*}

\noindent
 Thus, using the estimate \eqref{2e14a} in the above inequality, we have 	
\begin{equation*}\label{2e16}
\begin{split}
  \|U^{n}\| \, &\le C \, \Big( 1 + \|U^0\| \Big),
\end{split}
\end{equation*}
where $C = 2 E_{\alpha}\bigg(\frac{11}{4} \, t_{n}^{\alpha}\bigg) \, \max \bigg\{ 1, \,  \Big( \frac{11 \, T^{\alpha}}{4 \, \Gamma{(1+\alpha)}} + 1\Big) \, \max\limits_{1 \le j \le n} \|f^{j-\sigma}\| \bigg\}$.\\

\noindent
To derive an estimate for $\|\nabla U^{n}\|$, we rewrite \eqref{2e6.a} using the operator $\Delta_h$ as follows:
\begin{equation}\label{2e78a}
\begin{split}
\left( {D}^{\alpha}_{N}U^{n-\sigma}, v_h \right) \, - \, a\big(l(U^{n, \sigma})\big) \, (\Delta_h U^{n, \sigma}, v_h) \, =& \, (f^{n-\sigma}, v_h), \; \,  \forall v_h \in V_h. \\
\end{split}
\end{equation}

\noindent
Now we take $v_h=\, - \Delta_h U^{n, \sigma}$ in $\eqref{2e78a}$ and then proceed in a similar manner as above to obtain
\begin{equation*}\label{2e86}
\begin{split}
  \|\nabla U^{n}\| \, &\le C \, \Big( 1 + \|\nabla U^0\| \Big),
\end{split}
\end{equation*}
where $C = 2 \, E_{\alpha}\Big(\frac{11}{2} t_{n}^{\alpha}\Big) \, \max \left\lbrace 1, \, \frac{1}{\sqrt{m_1}} \Big( \frac{11 \, T^{\alpha}}{4 \, \Gamma{(1+\alpha)}} + 1\Big) \, \max\limits_{1 \le j \le n} \|f^{j-\sigma}\| \right\rbrace$ and the constant $m_1$ is as per Hypothesis H1 in the introduction. This completes the proof. \hfill $\square$

\section{Existence-Uniqueness of Fully-discrete Solution}

 \indent
 For proving the existence and the uniqueness of a fully-discrete solution, we shall need the following proposition, which is a consequence of the Brouwer's fixed point theorem \cite[P.237]{[vth]}.

\begin{prop}\label{p1}
	Let $S:X\rightarrow X$ be a continuous map on a finite dimensional Hilbert space X such that $(S(w),w)>0, \; \forall w \in X \; \mbox{with} \; \|w\|_X=\zeta$, $\zeta > 0$. Then there exists $v\in X$ such that $S(v)=0, \; \; \|v\|_X \leq \zeta.$
\end{prop}

We also define
\begin{equation}\label{2e120}
  \tau := \frac{2m_1 \, (a_{n,n} + b_{n,n-1})}{(1-\sigma) \, \{ L C R_1 (1 -2 \sigma) \}^2},
\end{equation}
where $R_1 = C (1 + \|\nabla U^0\|)$.\\

\indent
In next result, we prove the existence and uniqueness of a fully-discrete solution. The proof is almost similar to the proof of~\cite[Theorem 4.1]{[sp1]}. However, we prove it here for the sake of completeness.

\begin{theorem}\label{t1}
Let $U^0, \ldots, U^{n-1}$ be given. Then for all $1 \le n \le N$, there exists a unique solution $U^{n} \in V_h$ of \eqref{2e6}, provided $\max\limits_{\, 1 \le n \le N } \, \tau_{n} \le \tau$.
\end{theorem}

\noindent \textbf{Proof:} First we prove the existence. Rewriting \eqref{2e6.a} as follows:
\begin{eqnarray}\label{2e87}
  ( U^{n}, \, v_h) \!\!\! & + &\!\!\! \frac{1}{g_{n,n}} \, a\big(l(U^{n, \sigma})\big) (\nabla U^{n, \sigma}, \, \nabla v_h) - \frac{1}{g_{n,n}} \, \big(f^{n - \sigma}, \, v_h \big) \nonumber \\
  \!\!\!& - &\!\!\! \frac{1}{g_{n,n}} \, \Big(g_{n,1} U^0 + \sum_{k=2}^{n} (g_{n,k} - g_{n,k-1}) U^{k-1}, \, v_h \Big) = 0.
\end{eqnarray}

\noindent
Now, we define a map $G : V_h \longrightarrow V_h$ such that
\begin{eqnarray}\label{2e88}
  \big( G(X^{n}), v_h \big) = ( X^{n}, \! \! \!  &v_h& \! \! \! \! \!) - \frac{1}{g_{n,n}} \, \Big(g_{n,1} U^0 + \sum_{k=2}^{n} (g_{n,k} - g_{n,k-1}) U^{k-1}, \, v_h \Big) \qquad \nonumber\\
  & + & \! \! \! \frac{(1-\sigma)}{g_{n,n}} \, a\big(l(U^{n, \sigma})\big) (\nabla X^{n}, \, \nabla v_h) - \frac{1}{g_{n,n}} \, \big(f^{n - \sigma}, \, v_h \big) \nonumber \\
  & + & \! \! \! \frac{ \sigma}{g_{n,n}} \, a\big(l(U^{n, \sigma})\big) (\nabla U^{n-1}, \, \nabla v_h).
\end{eqnarray}
Then the map $G$ is continuous. By choosing $v_h = X^{n}$ in \eqref{2e88}, we get

\noindent
\begin{eqnarray}\label{2e91}
  \big( G(X^{n}), X^{n} \big) \ge   \! \! \! &\Big(& \! \! \!  \! \! \! \| X^{n} \| - \frac{g_{n,1}}{g_{n,n}} \| U^0 \| - \frac{1}{g_{n,n}} \, \sum_{k=2}^{n} (g_{n,k} - g_{n,k-1}) \, \| U^{k-1} \|  \nonumber\\
   & - & \! \! \! \frac{ \sigma \, m_1}{g_{n,n}} \, \| \nabla U^{n-1} \| - \frac{1}{g_{n,n}} \, \|f^{n - \sigma}\| \Big) \, \| X^{n} \|.
\end{eqnarray}

\noindent
Thus, $\big( G(X^{n}), X^{n} \big) > 0$, if
\begin{eqnarray}\label{2e93}
  \| X^{n} \| \, > \, \frac{g_{n,1}}{g_{n,n}} \| U^0 \| \! \! \! & + & \! \! \! \frac{1}{g_{n,n}} \, \sum_{k=2}^{n} (g_{n,k} - g_{n,k-1}) \, \| U^{k-1} \| \nonumber \\
  & + & \! \! \! \frac{ \sigma \, m_1}{g_{n,n}} \, \| \nabla U^{n-1} \| + \frac{1}{g_{n,n}} \, \|f^{n - \sigma}\|.
\end{eqnarray}

\noindent
Hence, by Proposition \ref{p1}, we can assure that \eqref{2e6} has a solution.\\

\noindent
Next we prove the uniqueness. Suppose, if possible, there exist two solutions of \eqref{2e6}, say $U^{n}_1$ and $U^{n}_2$ at time $t= t_n$. For simplicity, we denote $U^{n}_1$ by $U_1$, $U^{n}_2$ by $U_2$, $U^{n, \sigma}_1$ by $U^{\sigma}_1$ and $U^{n, \sigma}_2$ by $U^{\sigma}_2$, respectively. Let $U_1 - U_2 = \mathcal{Y}$.\\
From \eqref{2e87}, we can get
\begin{eqnarray}\label{2e95}
  ( U_1 - U_2, \, v_h) \!\!\!& + &\!\!\! \frac{(1 - \sigma)}{g_{n,n}} \, a\big(l(U^{\sigma}_1)\big) (\nabla U_1, \, \nabla v_h)  - \frac{(1 - \sigma)}{g_{n,n}} \, a\big(l(U^{\sigma}_2)\big) (\nabla U_2, \, \nabla v_h)  \nonumber \\
  \!\!\!& + &\!\!\! \frac{ \sigma}{g_{n,n}} \, \Big\{ a\big(l(U^{\sigma}_1)\big) - a\big(l(U^{\sigma}_2)\big) \Big\} \, (\nabla U^{n-1}, \nabla v_h) = 0.
\end{eqnarray}

\noindent
By subtracting $\frac{(1 - \sigma)}{g_{n,n}} \, a\big(l(U^{\sigma}_1)\big) \, (\nabla U_2, \nabla v_h)$ from both sides of \eqref{2e95} and then choosing $v_h = \mathcal{Y}$, we get
\begin{equation*}\label{2e97}
\begin{split}
( \mathcal{Y}, \, \mathcal{Y}) + \frac{(1 - \sigma)}{g_{n,n}} \, a\big(l(U^{\sigma}_1)\big) (\nabla \mathcal{Y}, \, \nabla \mathcal{Y})  =  & \, \frac{(1 - \sigma)}{g_{n,n}} \Big\{ a\big(l(U^{\sigma}_2)\big) - a\big(l(U^{\sigma}_1)\big) \Big\} (\nabla U_2, \, \nabla \mathcal{Y})  \qquad \qquad \\
& -  \frac{ \sigma}{g_{n,n}} \, \Big\{ a\big(l(U^{\sigma}_1)\big) - a\big(l(U^{\sigma}_2)\big) \Big\} (\nabla U^{n-1}, \nabla \mathcal{Y}).
\end{split}
\end{equation*}

\noindent
By using \eqref{2e102}, Hypotheses H1-H2, the Cauchy-Schwartz inequality, and the Young's inequality in the above equation, we have
\begin{equation*}\label{2e10}
\begin{split}
  \Big(1-\frac{(1 - \sigma) \{L C R_1 (1 - 2 \sigma)\}^2}{2 \, m_1 \, g_{n,n}} \Big) \, \| \mathcal{Y} \|^2  \, \le 0.\\
\end{split}
\end{equation*}

\noindent
Since $\max\limits_{\, 1 \le n \le N } \, \tau_{n} \le \tau$, it follows that $\Big(1-\frac{(1 - \sigma) \{L C R_1 (1 - 2 \sigma)\}^2}{2 \, m_1 \, g_{n,n}} \Big) > 0$. Hence, $U_1 = U_2.$
This completes the proof.  \hfill $\square$

\section{Error Estimates}\label{2err_est}
\indent
For deriving \emph{a priori} error estimates in Theorem \ref{2th5}, the following definition and lemma will be required.

 \begin{definition}\cite[P.8]{[vth]}
 	The Ritz-projection is a map  $R_h : H^1_0(\Omega) \rightarrow V_h$ such that
 	\begin{equation*}\label{2e62}
 	(\nabla w, \, \nabla v) \, = \, (\nabla R_h w, \, \nabla v), \quad \forall w \in  H^1_0(\Omega) \: \mbox{and}  \; \;  \forall v \in V_h.
 	\end{equation*}
 \end{definition}
 
 \indent
 Utilizing \eqref{2e60}, one can easily prove that $R_h $  satisfies $\|{\nabla R_h w}\| \, \le \, C$.
 
 \indent
 Since $V_h$ is the space of linear polynomials, it is well known \cite[Lemma 1.1]{[vth]} that
 \begin{equation}\label{2e61}
 \begin{split}
 \|{w-R_hw} \| \, + h \, \|{\nabla (w-R_hw)}\| \, \le& \, Ch^2 \, \|{\Delta w}\|, \quad \forall w \in H^2 \cap H^1_0.\\
 \end{split}
 \end{equation}

\begin{lemma}\label{2l3}
	\cite[P.9 (32)]{[hr12]} If $u \in C^2(0,T]$ and $\| \partial_t^{q} u\| \le C \, (1+t^{\alpha - q})$, for $q = 0,1,2,3$, then, for $n=1, \ldots, N,$ 	
	\begin{equation*}\label{2e123}
	\begin{split}
	t_{n-\sigma}^{\alpha} \, \|u^{n, \sigma} - u^{n-\sigma} \| \le C N^{-\min\{r \alpha, \, 2\}}.
	\end{split}
	\end{equation*}
\end{lemma}

\indent
Now, in the following theorem, we state and prove the most significant result of the present work. For this, we have to assume that the solution $u$ satisfies the following conditions.
\begin{eqnarray}\label{2e60}
&{^{c}_{0}{D}^{\alpha}_{t} u}& \! \! \! \in {L^{\infty}(0, T; {H^2(\Omega)})}, \qquad   u \in {L^{\infty}(0, T; {H^1_0(\Omega) \cap H^3(\Omega)})}, \nonumber \\
&\! \! \! \mbox{and}&  \, \| \partial_t^{q} u\|_3 \le C \, (1+t^{\alpha - q}), \quad \mbox{for} \;\, q = 0,1,2,3. 
\end{eqnarray}
\vspace{-20pt}

\begin{theorem}\label{2th5}
For $1 \le n \le N$, let $u^{n}$ and $U^{n}$ be the solutions of \eqref{2e1} and \eqref{2e6}, respectively. Then
\begin{eqnarray}
\label{2e20}
  \max_{1 \le n \le N} \|u^{n} - U^{n}\|& \le & C \, \big(h^2 + N^{-\min \left\lbrace 2, \, r \alpha \right\rbrace }\big), \\
\label{2e21}
  \max_{1 \le n \le N} \|\nabla (u^{n} - U^{n}) \| & \le & C \, \big(h + N^{-\min \left\lbrace 2, \, r \alpha \right\rbrace }\big).
\end{eqnarray}
\end{theorem}

\noindent \textbf{Proof.}
For the proof, we split the error ($u^n - U^n$) as follows:
\begin{equation*}\label{2e61n}
\begin{split}
u^n - U^n = u^n - R_hu^n + R_hu^n - U^n = \rho ^n + \theta ^n,
\end{split}
\end{equation*}
where \  $ \rho ^n := u^n - R_hu^n$ \ and \ $\theta ^n := R_hu^n - U^n.$
As the estimate for $ \rho ^n$ can be given by \eqref{2e61}, it is sufficient to find estimate for $\theta ^n$ for proving the theorem. \\

\noindent
For any $v_h \in V_h$, the estimate for $\theta ^n$ is given by
\begin{eqnarray}\label{2e22}
  \big( {D}^{\alpha}_{N} \theta^{n-\sigma}, \! \! \! \! &v_h& \! \! \! \! \! \! \! \big) \, + \, a\big(l(U^{n, \sigma})\big) \, (\nabla \theta^{n, \sigma}, \nabla v_h) \nonumber\\
  & = & \! \! \! \! \! \big( {D}^{\alpha}_{N} R_hu^{n-\sigma} - \, ^{c}_{0}{D}^{\alpha}_{t_{n-\sigma}}u , v_h \big) - a\big(l(u^{n-\sigma})\big) (\Delta u^{n, \sigma} - \Delta u^{n-\sigma}, v_h) \nonumber \\
  & \quad - & \! \! \! \!  \big\{ a\big(l(U^{n, \sigma})\big) - a\big(l(u^{n-\sigma})\big) \big\} (\Delta u^{n, \sigma}, v_h).
\end{eqnarray}

\noindent
Now, we choose $v_h = \theta^{n, \sigma}$ in \eqref{2e22} to get
\begin{eqnarray}\label{2e23}
  \big( {D}^{\alpha}_{N} \theta^{n-\sigma}, \theta^{n, \sigma} \big) + a\big(l(U^{n, \sigma})\big) \, (\nabla \theta^{n, \sigma}, \nabla \theta^{n, \sigma}) \! \! \! &=& \! \! \! \big( {D}^{\alpha}_{N} R_hu^{n-\sigma} - \, ^{c}_{0}{D}^{\alpha}_{t_{n-\sigma}}u , \theta^{n, \sigma} \big) \nonumber\\
   &-& \! \! \! \big\{ a\big(l(U^{n, \sigma})\big) - a\big(l(u^{n-\sigma})\big) \big\} (\Delta u^{n, \sigma}, \theta^{n, \sigma})  \nonumber \\
  &-& \! \! \! a\big(l(u^{n-\sigma})\big) (\Delta u^{n, \sigma} - \Delta u^{n-\sigma}, \theta^{n, \sigma}).
\end{eqnarray}

\noindent
An application of the bound of $a$, the Cauchy-Schwartz inequality and the Triangle inequality in \eqref{2e23} gives	
\begin{eqnarray}\label{2e24}
  \big( {D}^{\alpha}_{N} \theta^{n-\sigma}, \theta^{n, \sigma} \big)  +  m_1 \|\nabla \theta^{n, \sigma}\|^2 & \le & \|{D}^{\alpha}_{N} R_hu^{n-\sigma} - ^{c}_{0}{D}^{\alpha}_{t_{n-\sigma}}R_hu \| \, \|\theta^{n, \sigma}\| \nonumber\\
  &  & + \|^{c}_{0}{D}^{\alpha}_{t_{n-\sigma}}R_hu - \, ^{c}_{0}{D}^{\alpha}_{t_{n-\sigma}}u\| \|\theta^{n, \sigma}\| \nonumber\\
  &  & + m_2 \, \|\Delta u^{n, \sigma} - \Delta u^{n-\sigma} \| \|\theta^{n, \sigma}\| \nonumber \\
  &  & + \big|a\big(l(U^{n, \sigma})\big) - a\big(l(u^{n-\sigma})\big)\big| \| \Delta u^{n, \sigma}\| \|\theta^{n, \sigma}\|.
\end{eqnarray}

\noindent
Since function $a$ is Lipschitz continuous, we have
\begin{eqnarray}\label{2e28}
  \big|a\big(l(U^{n, \sigma})\big) - a\big(l(u^{n-\sigma})\big)\big| &\le& L \, C \, \|U^{n, \sigma} - u^{n-\sigma}\| \nonumber \\
  &\le& L \, C \, \big\{ \| U^{n, \sigma} - u^{n, \sigma}\| + \|u^{n, \sigma} - u^{n-\sigma}\| \big\}  \nonumber \\
  &\le& L \, C \, \big\{ \| \rho^{n, \sigma}\| + \|\theta^{n, \sigma}\| + \|u^{n, \sigma} - u^{n-\sigma}\| \big\}.
\end{eqnarray}

\noindent
Also, it follows from \eqref{2e60} and \eqref{2e61} that
\begin{equation}\label{2e31}
\begin{split}
  \|^{c}_{0}{D}^{\alpha}_{t_{n-\sigma}}R_hu - \, ^{c}_{0}{D}^{\alpha}_{t_{n-\sigma}}u\| \; \le \; C \, h^2 \, \| \Delta \, ^{c}_{0}{D}^{\alpha}_{t_{n-\sigma}}u \| \; \le \; C \, h^2.
\end{split}
\end{equation}

\noindent
From Lemma \ref{2l1}, \eqref{2e61}, \eqref{2e24}, \eqref{2e28} and \eqref{2e31}, we have
\begin{eqnarray}\label{2e32}
  \! \! \! \! \! \! {D}^{\alpha}_{N} \|\theta^{n-\sigma}\|^2 \le C_1 \, \|\theta^{n, \sigma}\|^2 +  C_2 \Big\{ \|{D}^{\alpha}_{N} R_hu^{n-\sigma} \! \! \! \!  &-& \! \! \! ^{c}_{0}{D}^{\alpha}_{t_{n-\sigma}}R_hu \| + \|u^{n, \sigma} - u^{n-\sigma}\| \nonumber \\
  & + & \! \! \!  \|\Delta u^{n, \sigma} - \Delta u^{n-\sigma}\| +h^2 \Big\} \|\theta^{n, \sigma}\|, \quad
\end{eqnarray}
where $C_1$ is dependent on $C$, $L$ and $C_2$ is dependent on $C$, $L$, $m_2$.\\

\noindent
An application of Lemma \ref{2l5} into \eqref{2e32} leads to
\begin{eqnarray}\label{2e107}
   \! \!  \!  \|\theta^{n}\| \le 2 \, E_{\alpha}\Big(\frac{11}{2} C_1 \! \! \! \! &t_{n}^{\alpha}& \! \! \! \! \! \Big) \Big( \|\theta^0\| + C_2 \max_{1 \le k \le n} \sum_{j=1}^{k} p^{(k)}_{k-j} \Big\{ h^2 + \|u^{j, \sigma} - u^{j- \sigma}\| \nonumber\\
   &+& \! \! \! \!  \! \! \|{D}^{\alpha}_{N} R_hu^{j- \sigma} - \, ^{c}_{0}{D}^{\alpha}_{t_{j- \sigma}}R_hu \| + \|\Delta u^{j, \sigma} - \Delta u^{j- \sigma}\| \Big\} \Big). \quad
\end{eqnarray}

\noindent
Now, from Lemmas~ \ref{2l6} and  \ref{2l3}, we have
\begin{eqnarray}\label{2e107a}
  \! \! \! \!&\|&\! \! \! \! \! {D}^{\alpha}_{N} R_hu^{j- \sigma} - \, ^{c}_{0}{D}^{\alpha}_{t_{j- \sigma}}R_hu \| + \|u^{j, \sigma} - u^{j- \sigma}\| + \|\Delta u^{j, \sigma} - \Delta u^{j- \sigma}\| \nonumber\\
 &\le&\! \! \! C \, t^{-\alpha}_{j-\sigma} \, \big(N^{-\min \{3- \alpha, \, r \alpha\}} + N^{- \min \{2, \, r \alpha\}}\big) \nonumber\\
 &\le&\! \! \! C \, (2N)^{r \alpha} \, N^{- \min \{2, \, r \alpha\}} \, j^{r (l_N - \alpha)},
\end{eqnarray}
%
where the fact that $t^{-\alpha}_{j-\sigma} \le j^{-r \alpha} (2N)^{r \alpha} \le (2N)^{r \alpha} j^{r (l_N - \alpha)}$ is used with $l_N = \frac{1}{\ln N}$.\\

\noindent
Using \eqref{2e107a} in \eqref{2e107}, we can obtain
\begin{eqnarray}\label{2e108}
  \|\theta^{n}\| \le 2 \, E_{\alpha}\Big(\frac{11}{2} \! \! \! \!&C_1&\! \! \! \!  t_{n}^{\alpha}\Big) \Big( \|\theta^0\| + C_2 h^2 \max_{1 \le k \le n} \sum_{j=1}^{k} p^{(k)}_{k-j} \nonumber\\
  & +&\! \! \! \! C \, (2N)^{r \alpha} N^{- \min \{2, \, r \alpha\}} \max_{1 \le k \le n} \sum_{j=1}^{k} p^{(k)}_{k-j} \, j^{r (l_N - \alpha)} \Big).
\end{eqnarray}

\noindent
Now, using Lemma \ref{2l7} with $\gamma = l_N$, we can have
\begin{equation}\label{2e108a}
\begin{split}
  \sum_{j=1}^{k} p^{(k)}_{k-j} \, j^{r (l_N - \alpha)} \, \le& \, \frac{11 \, \Gamma(1 + l_n - \alpha)}{4 \, \Gamma{(1+l_N)}} \; T^{\alpha} \, N^{r(l_N-\alpha)}.
\end{split}
\end{equation}

\noindent
Using \eqref{2e14a} and \eqref{2e108a} in \eqref{2e108}, we get
\begin{eqnarray}\label{2e111}
\|\theta^{n}\| \le 2 \, E_{\alpha}\Big(\frac{11}{2} \! \! \! \!&C_1&\! \! \! \! \, t_{n}^{\alpha}\Big) \Big( \|\theta^0\| + \frac{11 \, C_2 \, T^{\alpha}}{4 \, \Gamma{(1+\alpha)}} \; h^2 \nonumber \\
& + &\! \! \! \! C \, (2N)^{r \alpha} \frac{11 \, \Gamma(1 + l_n - \alpha)}{4 \, \Gamma{(1+l_N)}} \; T^{\alpha} \, N^{r(l_N-\alpha)} \,N^{- \min \{2, \, r \alpha\}} \Big).
\end{eqnarray}

\noindent
Choosing $U^0  = R_hu_0$, we get $\|\theta^0\| = 0$. Therefore, from \eqref{2e111}, we have,
\begin{equation}\label{2e112}
\begin{split}
  \|\theta^{n}\| &\le C \, \big( h^2 + N^{- \min \{2, \, r \alpha\}} \big),\\
\end{split}
\end{equation}
where $C = 2 \, E_{\alpha}\Big(\frac{11}{2} C_1 \, t_{n}^{\alpha}\Big) \, \max \Big\{ \frac{11 \, C_2 \, T^{\alpha}}{4 \, \Gamma{(1+\alpha)}}, \, 2^{r \alpha} \, C  \frac{11 \, \Gamma(1 + l_n - \alpha)}{4 \, \Gamma{(1+l_N)}} \; N^{r l_N} \, T^{\alpha} \Big\}$.\\

\noindent
Finally, an application of the Triangle inequality and the estimates \eqref{2e61} and \eqref{2e112} give \eqref{2e20}.\\

\noindent
Now, to derive the error estimate in $H^1_0$-norm, we use the definition of discrete Laplacian $\Delta_h$ and rewrite \eqref{2e22} as follows:
\begin{eqnarray}\label{2e33}
\big( {D}^{\alpha}_{N} \theta^{n-\sigma}, \! \! \! \! &v_h& \! \! \! \! \! \! \! \big) \, - \, a\big(l(U^{n, \sigma})\big) \, (\Delta_h \theta^{n, \sigma}, v_h) \nonumber\\
& = & \! \! \! \! \! \big( {D}^{\alpha}_{N} R_hu^{n-\sigma} - \, ^{c}_{0}{D}^{\alpha}_{t_{n-\sigma}}u , v_h \big) - a\big(l(u^{n-\sigma})\big) (\Delta u^{n, \sigma} - \Delta u^{n-\sigma}, v_h) \nonumber \\
& \quad - & \! \! \! \!  \big\{ a\big(l(U^{n, \sigma})\big) - a\big(l(u^{n-\sigma})\big) \big\} (\Delta u^{n, \sigma}, v_h).
\end{eqnarray}
\noindent
We take $v_h = - \Delta_h \theta^{n, \sigma}$ in \eqref{2e33} and then use the similar arguments as above together with the Poincar$\acute{e}$ inequality to obtain
\begin{equation*}\label{2e118}
\begin{split}
  \|\nabla  u^{n} - \nabla  U^{n}\|  \le& \, C \, \big( h + N^{- \min \{2, \, r \alpha\}} \big).\\
\end{split}
\end{equation*}
This completes the proof. \hfill $\square$

\begin{corollary}
For $1 \le n \le N$, let $u^{n}$ and $U^{n}$ be the solutions of \eqref{2e1} and \eqref{2e6}, respectively. Then
\begin{equation*}
  \max_{1 \le n \le N} \|{u^{n} - U^{n}}\|_1 \leq C \, \big( h + {N}^{-\min \left\lbrace 2, \, r \alpha \right\rbrace }  \big).
\end{equation*}
\end{corollary}
\noindent \textbf{Proof.} This proof follows from the Poincar$\acute{e}$ inequality and Theorem \ref{2th5}.  \hfill $\square$\\

\begin{remark}
      Above analysis can be easily extended for the case when $f=f(u)$ is Lipschitz continuous, except the $H_0^1$-norm error analysis. More precisely, in case of $f=f(u)$, we get a term $\| \nabla \big( f(u^{n - \sigma}) - f(U^{n, \sigma})\big) \|$ which creates problem in getting $H_0^1$-norm error estimate.
\end{remark}

\section{Numerical Experiments}
\indent In this section, for confirmation of our theoretical estimates, we perform two numerical experiments with known exact solution. In both numerical experiments, the time interval is taken to be $[0,1]$. In order to achieve optimal order of convergence rate, we choose mesh grading parameter $r=\frac{2}{\alpha}.$ We consider the tolerance $\epsilon = 10^{-12}$ for stopping the iterations in Newton's method. Spatial domain $\Omega$ is partitioned with $\big(M_s+1\big)$ node points in each direction and the time interval is partitioned into $N$ number of sub-intervals. We set $M_s=N$ and calculate error for different values of $N$ to get the rate of convergence in time direction. Similarly, we set $N=M_s$ and calculate error for different values of $M_s$ to conclude the rate of convergence in spatial direction. \\

\begin{example}\label{2E2}
	For this example, we consider equation \eqref{2e1} with a spatial domain $\Omega = (0, \pi)$ and $a(w)=3+\sin w$. Then we can choose $f(x,t)$ such that the analytical solution of equation \eqref{2e1} is $u(x,t)=(t^3+t^{\alpha}) \sin x$. 	 	
\end{example}

\indent 
Table \ref{2T5} shows the order of convergence in the temporal direction in $L^{\infty}(0,T;L^2(\Omega))$ norm on graded mesh for $\alpha = 0.4,$ $\alpha = 0.6$ and $\alpha = 0.8$. From this table, it can be seen that using graded mesh, we can achieve $O(N^{-2})$ convergence rate in the temporal direction. The order of convergence in the spatial direction in $L^2$-norm and $H^1_0$-norm for $\alpha = 0.4,$ $\alpha = 0.6$ and $\alpha = 0.8$ are given in Tables \ref{2T6} and \ref{2T7}, respectively. From these tables, one can observe that the numerical convergence results are in accordance with our theoretical convergence estimates. 
\newpage
\begin{center}
	\begin{tiny}
		\begin{table}[h!]
			\begin{center}
				\begin{tabular}{|c|c|c|c|c|c|c|}
					\hline
					\multirow{2}{*}{\large $N$} & \multicolumn{2}{c|}{$\alpha = 0.4$} & \multicolumn{2}{c|}{$\alpha = 0.6$} & \multicolumn{2}{c|}{$\alpha = 0.8$}\\
					\cline{2-7}
					& Error & OC & Error & OC & Error & OC \\
					\hline
					$2^6$ & 1.22E-03 & 1.9748 & 1.04E-03 & 2.0270 & 1.01E-03 & 2.0136 \\
					\hline
					$2^7$ & 3.11E-04 & 2.0069 & 2.54E-04 & 2.0239 & 2.49E-04 & 2.0143 \\
					\hline
					$2^8$ & 7.75E-05 & 2.0073 & 6.26E-05 & 2.0149 & 6.17E-05 & 2.0138 \\
					\hline
					$2^9$ & 1.93E-05 & 2.0068 & 1.55E-05 & 2.0125 & 1.53E-05 & 2.0128 \\ 
					\hline
					$2^{10}$ & 4.80E-06 & - & 3.84E-06 & - & 3.79E-06 & - \\
					\hline 	
				\end{tabular}
			\end{center}
			\caption {\emph {Error and order of convergence  in $L^{\infty}(0,T;L^2(\Omega))$ norm in temporal direction on graded mesh for Example \ref{2E2}.}}
			\label{2T5}
		\end{table}
	\end{tiny}
\end{center}
\begin{center}
	\begin{tiny}
		\begin{table}[h!]
			\begin{center}	
				\renewcommand{\arraystretch}{1.2}
				\begin{tabular}{|c|c|c|c|c|c|c|}
					\hline
					\multirow{2}{*}{\large $M_s$} &    \multicolumn{2}{c|}{$\alpha = 0.4$} & \multicolumn{2}{c|}{$\alpha = 0.6$} & \multicolumn{2}{c|}{$\alpha = 0.8$} \\
					\cline{2-7}
					& Error & OC  & Error & OC & Error & OC  \\
					\hline
					$2^6$ & 7.97E-04 & 2.0785 & 1.04E-03 & 	2.0270 & 1.01E-03 & 2.0136 \\
					\hline
					$2^7$ & 1.89E-04 & 2.0778 & 2.54E-04 & 2.0274 & 2.49E-04 & 2.0143  \\
					\hline
					$2^8$ & 4.47E-05 & 2.0629 & 6.24E-05 & 2.0243 & 6.17E-05 & 2.0139 \\
					\hline
					$2^9$ & 1.07E-05 & 2.0472 & 1.53E-05 & 2.0203 & 1.53E-05 & 2.0128 \\ 
					\hline
					$2^{10}$ & 2.59E-06 & - & 3.78E-06 & - & 3.79E-06 & - \\
					\hline 
				\end{tabular}
				\caption {\emph{Error and order of convergence in $L^2$-norm in space for Example \ref{2E2}.}}
				\label{2T6}
			\end{center}
		\end{table}
	\end{tiny}
\end{center}
\begin{center}
	\begin{tiny}
		\begin{table}[h!]
			\begin{center}
				\renewcommand{\arraystretch}{1.2}
				\begin{tabular}{|c|c|c|c|c|c|c|}
					\hline
					\multirow{2}{*}{\large $M_s$} &    \multicolumn{2}{c|}{$\alpha = 0.4$} & \multicolumn{2}{c|}{$\alpha = 0.6$} & \multicolumn{2}{c|}{$\alpha = 0.8$} \\
					\cline{2-7}
					&  Error & OC &  Error & OC &  Error & OC  \\
					\hline
					$2^6$ & 3.55E-02 & 0.9999  & 3.55E-02 & 1.0000 & 3.55E-02 &  1.0000 \\
					\hline
					$2^7$ & 1.78E-02 & 0.9999 & 1.78E-02 & 1.0000 & 1.78E-02 & 1.0000  \\
					\hline
					$2^8$ & 8.88E-03 & 0.9999 & 8.88E-03 & 1.0000 & 8.88E-03 & 1.0000  \\
					\hline
					$2^9$ & 4.44E-03 & 0.9999 & 4.44E-03 & 1.0000 & 4.44E-03 & 1.0000 \\  
					\hline
					$2^{10}$ & 2.22E-03 & - & 2.22E-03 & - & 2.22E-03 & - \\
					\hline 					
				\end{tabular}
				\caption {\emph {Error and order of convergence  in $H^1_0$-norm in space for Example \ref{2E2}.}}
				\label{2T7}
			\end{center}
		\end{table}
	\end{tiny}
\end{center}

\begin{example}\label{2E3}
	In this second example, consider equation \eqref{2e1} with the spatial domain $\Omega = (0,1) \times (0,1)$ and $a(w)=3+\sin w$. Then we can choose $f(x, y, t)$ in such a way that the analytical solution of the equation \eqref{2e1} is $u(x,y,t)=(t^3+t^{\alpha}) (x-x^2)(y-y^2)$.
\end{example}
\newpage
\indent 
In Table \ref{2T8}, we have given the order of convergence in the temporal direction in $L^{\infty}(0,T;L^2(\Omega))$ norm on graded mesh for $\alpha = 0.5,$ $\alpha = 0.7$ and $\alpha = 0.9$. It can be seen from this table that using graded mesh, we obtain optimal convergence rate in time direction, as predicted by Theorem \ref{2th5}.   Tables \ref{2T9} and \ref{2T10} confirm the theoretical convergence estimates in $L^2$-norm \& $H^1_0$-norm in spatial direction for $\alpha = 0.5,$ $\alpha = 0.7$ and $\alpha = 0.9$. 
\begin{center}
	\begin{tiny}	
		\begin{table}[h!]
			\begin{center}
				\begin{tabular}{|c|c|c|c|c|c|c|}
					\hline
					\multirow{2}{*}{\large $N$} & \multicolumn{2}{c|}{$\alpha = 0.5$} & \multicolumn{2}{c|}{$\alpha = 0.7$} & \multicolumn{2}{c|}{$\alpha = 0.9$} \\
					\cline{2-7}
					& Error & OC & Error & OC & Error & OC  \\
					\hline
					$2^3$ & 5.59E-03 & 1.8107 & 4.94E-03 & 1.9108 & 4.40E-03 & 1.9553 \\
					\hline
					$2^4$ & 1.59E-03 & 1.9193 & 1.31E-03 & 1.9680 & 1.14E-03 & 1.9872 \\
					\hline
					$2^5$ & 4.21E-04 & 1.9644 & 3.36E-04 & 1.9875 & 2.86E-04 & 1.9960 \\
					\hline
					$2^6$ & 1.08E-04 & 1.9836 & 8.46E-05 & 1.9947 & 7.18E-05 & 1.9987 \\
					\hline
					$2^7$ & 2.73E-05 & - & 2.12E-05 & - & 1.80E-05 & - \\
					\hline															
				\end{tabular}
			\end{center}
			\caption {\emph{Error and order of convergence in $L^{\infty}(0,T;L^2(\Omega))$ norm in time on graded mesh for Example \ref{2E3}.}}
			\label{2T8}
		\end{table}
	\end{tiny}
\end{center}
\begin{center}
	\begin{tiny}	
		\begin{table}[h!]
			\renewcommand{\arraystretch}{1.2}
			\begin{center}
				\begin{tabular}{|c|c|c|c|c|c|c|}
					\hline
					\multirow{2}{*}{\large $M_s$} & \multicolumn{2}{c|}{$\alpha = 0.5$} & \multicolumn{2}{c|}{$\alpha = 0.7$} & \multicolumn{2}{c|}{$\alpha = 0.9$}\\
					\cline{2-7}
					& Error & OC & Error & OC & Error & OC  \\
					\hline
					$2^3$ & 5.59E-03 & 1.8107 & 4.94E-03 & 1.9108 & 4.40E-03 & 1.9553 \\
					\hline
					$2^4$ & 1.59E-03 & 1.9193 & 1.31E-03 & 1.9680 & 1.14E-03 & 1.9872 \\
					\hline
					$2^5$ & 4.21E-04 & 1.9644 & 3.36E-04 & 1.9875 & 2.86E-04 & 1.9960 \\
					\hline
					$2^6$ & 1.08E-04 & 1.9836 & 8.46E-05 & 1.9947 & 7.18E-05 & 1.9987 \\
					\hline
					$2^7$ & 2.73E-05 & - & 2.12E-05 & - & 1.80E-05 & - \\
					\hline 															
				\end{tabular}
			\end{center}
			\caption {\emph {Error and order of convergence in $L^2$-norm in space for Example \ref{2E3}.}}
			\label{2T9}
		\end{table}
	\end{tiny}
\end{center}
\newpage
\begin{center}
	\begin{tiny}	
		\begin{table}[h!]
			\renewcommand{\arraystretch}{1.2}
			\begin{center}
				\begin{tabular}{|c|c|c|c|c|c|c|}
					\hline
					\multirow{2}{*}{\large $M_s$} & \multicolumn{2}{c|}{$\alpha = 0.5$} & \multicolumn{2}{c|}{$\alpha = 0.7$} & \multicolumn{2}{c|}{$\alpha = 0.9$} \\
					\cline{2-7}
					& Error & OC & Error & OC & Error & OC \\
					\hline
					$2^3$ & 4.67E-02 & 1.0359 & 4.59E-02 & 1.0232 & 4.53E-02 & 1.0120  \\
					\hline
					$2^4$ & 2.28E-02 & 1.0171 & 2.26E-02 & 1.0080 & 2.25E-02 & 1.0035 \\
					\hline
					$2^5$ & 1.13E-02 & 1.0055 & 1.12E-02 & 1.0022 & 1.12E-02 & 1.0009  \\
					\hline
					$2^6$ & 5.61E-03 & 1.0015 & 5.60E-03 & 1.0006 & 5.60E-03 & 1.0002 \\  
					\hline 
					$2^7$ & 2.80E-03 & - & 2.80E-03 & - & 2.80E-03 & - \\  
					\hline 					       		    
				\end{tabular}
			\end{center}
			\caption {\emph {Error and order of convergence in $H^1_0$-norm in space for Example \ref{2E3}.}}
			\label{2T10}
		\end{table}
	\end{tiny}
\end{center}

\section{Conclusions}
In this work, we have proposed $L2$-$1_{\sigma}$ Galerkin finite element method for solving time-fractional PDE with nonlocal diffusion term. Under certain regularity assumptions on the exact
solution, we have derived $\alpha$-robust \emph{a priori} error estimates in $L^2$-norm and $H^1$-norm. We also have confirmed our theoretical estimates by numerical experiments. 
The derivation of regularity assumptions (in Section \ref{2err_est}) on the exact solution $u$ is under the investigation and it will be reported in our future work.  \\


\noindent \textbf{Declarations:}\\
\textbf{Conflict of interest-} The author declares no competing interests.

\end{document}